\documentstyle[bezier,pstricks,pst-node,numreferences]{kluwer}

\def\komment#1{}

\font\yswab=yswab scaled 1200
\font\twlsf=cmssi12
\font\tensf=cmssi10
\font\eigsf=cmssqi8
\font\bigti=cmss12
\font\twlBbb=msbm10 scaled 1200

\font\ninBbb=msbm10 scaled 900

\font\sevBbb=msbm10 scaled 700
\newfam\Bbbfam
\def\Bbb{\fam\Bbbfam \twlBbb}
\textfont\Bbbfam\twlBbb \scriptfont\Bbbfam\ninBbb \scriptscriptfont\Bbbfam\sevBbb

\newcommand{\altn}{\mbox{{\yswab n}}{}}
\newcommand{\altg}{\mbox{{\yswab g}}{}}
\komment{
\def\reell{{\mathchoice{{I\!\!R}}{{I\!\!R}}{{\scriptstyle I\!\!R}}%
{{\scriptscriptstyle I\!\!R}}}}
\newfam\ssffam

\textfont\ssffam\twlsf \scriptfont\ssffam\tensf \scriptscriptfont\ssffam\eigsf
\def\Z{{\ssf Z\!\!Z}}
\def\natur{{\mathchoice{{I\!\!N}}{{I\!\!N}}{{\scriptstyle I\!\!N}}%
{{\scriptscriptstyle I\!\!N}}}}
}
\def\Z{{\Bbb Z}}
\def\reell{{\Bbb R}}
\def\natur{{\Bbb N}}
\def\ga{{\gamma}}

\def\al{{\alpha}}

\def\ep{{\epsilon}}

\def\om{{\omega}}

\def\De{{\Delta}}
\def\Ga{{\Gamma}}

\def\LL{{\cal L}}

\def\HH{{\cal H}}
\def\PP{{\cal P}}

\def\CC{{\cal C}}
\def\RR{{\cal R}}
 
\def\ti{\tilde}
\def\tiM{\widetilde{\hbox{$M$}}}
\def\tiMinf{{\tiM}_{\!\infty}}

\def\vp{\vec{p}}
\def\vq{\vec{q}}

\def\SGA{{S_\Ga\dot{\cup}S_\Ga^{-1}}}
\def\word#1{{(#1_1, #1_2, \dots)}}
\def\wordl#1#2{{(#1_1, #1_2, \dots,#1_#2)}}
\def\la{{\langle}}
\def\ra{{\rangle}}
\def\hei#1{{\HH_{#1}}}
\def\ohne{\setminus}

\def\leerklam{{\mathchoice{{\la \, .\, , .\, \ra}}{{\la \, .\, , .\, \ra}}{{\la \, .\, ,\, .\, \ra}}{{\la \, .\, , \, .\, \ra}}}}
\def\zvec#1#2{\pmatrix{\hfill #1\cr \hfill #2}}
\def\empt{{\emptyset}}
\def\ddal{{\res{\partial\over \partial \al}{\al=0}}}

\def\kernel{\mathop{\rm ker}}

\def\dim{\mathop{\rm dim}}
\def\rank{\mathop{\rm rank}}
\def\injrad{\mathop{\rm inj\,rad}}
\def\area{\mathop{\rm Area}}
\def\length{\LL}
\def\bigtimes{\mathop{\hbox{\bigti X}}\limits}
\def\Konst{{\mathchoice{{K}}{{K}}{{\scriptstyle K}}{{\scriptscriptstyle K}}}}
\def\vf#1{{{\cal V}_{#1}}}

\def\qed{{\leavevmode\unskip\nobreak\hfil\penalty 50\hskip 1em%
  \hbox{}\nobreak\hfil\lower 1pt\hbox{$\Box$\kern-.5pt}\parfillskip 0pt
  \finalhyphendemerits 0\par}}
\def\me#1#2 {\ifcat,#2{\em #1}\else {\em #1\/}\fi #2\ }
\def\res#1#2{{#1}\lower .11ex\hbox{$|$}\lower .644ex\hbox{$\scriptstyle #2$}}

\def\fatlines{\linethickness{1pt}}
\def\smalllines{\linethickness{0.4pt}}

\newtheorem{theorem}{THEOREM}[section]
\newtheorem{definition}[theorem]{DEFINITION}

\newtheorem{proposition}[theorem]{PROPOSITION}
\newtheorem{lemma}[theorem]{LEMMA}
\newtheorem{corollary}[theorem]{COROLLARY}
\newtheorem{openproblem}[theorem]{OPEN PROBLEM}
\newtheorem{remark}[theorem]{Remark}
\newtheorem{supplement}[theorem]{SUPPLEMENT}

\def\spacebetweentheoremandproof{\vskip3mm}
\def\head#1{\par\smallskip\noindent{\it #1\ }}
\def\remark{{\head{Remark.}}}

\newlength{\dotstrength}
\newlength{\curvestrength}
\newlength{\framestrength}
\newlength{\pointerstrength}
\newlength{\fibrestrength}
\newlength{\anglestrength}
\newlength{\angledotstrength}
\begin{opening}
\title{Minimal Geodesics and Nilpotent Fundamental Groups}
\author{Bernd \surname{Ammann}}
\institute{Mathematisches Institut der Universit\"at Freiburg\\
D-79104 Freiburg\\Germany\\
E-Mail: ammann@arcade.mathematik.uni-freiburg.de}
\date{March 1996}
\end{opening}

\begin{document}
\maketitle
\begin{abstract}
Hedlund \cite{hedlund} constructed Riemannian metrics on $n$-tori, $n \geq 3$
for which minimal geodesics are very rare. In this paper we construct 
similar examples for every nilpotent fundamental group.
These examples show that Bangert's existence results of minimal geodesics 
\cite{bangertmin} are optimal for nilpotent fundamental 
groups.
\end{abstract}

\classification{AMS Subject Classification (1991): Primary 53C22, Secondary 22E25, 20F18}


\section{Introduction}

\begin{definition}
Let $(M,\leerklam)$ be a Riemannian manifold. A non-constant geodesic
$c\colon \reell \to M$ is called \me{minimal} 
if it
satisfies for all $t_1<t_2$:
   $$\LL\left(\res{c}{[t_1,t_2]}\right) \leq 
     \LL\left(\ga\right)$$
for all curves $\gamma$ homotopic to
$\res{c}{[t_1,t_2]}$ with fixed endpoints.
\end{definition}

Suppose we fix $\pi_1(M)$. Then there are several known results that guarantee
the existence of minimal geodesics. 

\komment{
Existence and properties of minimal geodesics on compact Riemannian manifolds 
$(M,\leerklam)$ depend on the topology of $M$, in particular on $\pi_1(M)$, 
and on the Riemannian metric.
The simplest result relating $\pi_1(M)$ to the existence of minimal geodesics
is that $(M,\leerklam)$ carries a minimal geodesic if and only if 
$\pi_1(M)$ is infinite.
}

The simplest one is that $(M,\leerklam)$ 
carries a minimal geodesic if and only if 
$\pi_1(M)$ is infinite.

For some classes of differentiable manifolds certain existence properties
of minimal 
geodesics do not depend on the choice of Riemannian metric:
the bestknown cases are compact manifolds $M$ with hyperbolic 
fundamental groups. Here one can compactify the universal cover $\tiM$
of $M$ by a ``boundary at infinity'' $\tiMinf$. For every
Riemannian metric on $M$ the lift of a minimal geodesic to $\tiM$
converges for $t\to \pm \infty$ to two different points on $\tiMinf$
and, conversely, for each pair of different points on $\tiMinf$
there exists such a minimal geodesic 
(\cite{morse},\cite{busemann},\cite{klingenberg},\cite{gromov3}7.5).

The situation is similar on the $2$-torus $T^2= {\reell^2}/{\Z^2}$
where for every straight line in $\reell^2$ and every 
Riemannian metric on $T^2$ one finds a minimal geodesic whose lift stays at 
finite distance from the straight line 
(\cite{hedlund},\linebreak[0]\cite{biapol},\linebreak[0]\cite{bangertdyna}).

Surprisingly, the situation is completely different for an $n$-torus 
$T^n={\reell^ n}/{\Z^ n}$ if $n\geq 3$.
Here existence properties of minimal geodesics 
depend very much on the choice of the 
Riemannian metric:
for flat metrics every geodesic is minimal (and lifts to a straight line in 
$\reell^n$). On the other hand, there are the Hedlund metrics \cite{hedlund} 
on $T^n$, discussed in $\cite{bangertmin}$, where one has only $n$ 
periodic minimal ones.
So in these Hedlund examples minimal geodesics are very rare.
Using the language of dynamical systems one would say that the set of 
unit tangent vectors to minimal geodesics 
consists of $2n$ periodic orbits of the 
geodesic flow and (countably many) heteroclinic and homoclinic connections
between them.

These Hedlund metrics contrast to a theorem of V.~Bangert 
(\cite{bangertmin},\cite{bangerticm}). He proves that the number 
of ``directions'' of minimal geodesics on an arbitrary Riemannian manifold
$(M,\leerklam)$ is at least the first Betti-number 
$b_1:=\rank \pi_1/[\pi_1,\pi_1]$.
This bound is 
optimal for the Hedlund metric on $T^n,n\geq 3$.
\komment{
But in the case that the fundamental group is very far from being abelian,
it seems that it should be possible to 
improve this bound. If e.g.\ $\pi_1(M)=Sl(3,\Z)$ then $M$ has
a minimal geodesic, but $b_1=0$. 
}

In this paper we  
will construct Riemannian manifolds with 
only $b_1$ different ``directions'' of minimal geodesics for arbitrary
nilpotent fundamental groups. Therefore Bangert's bound is optimal 
for arbitrary nilpotent groups.

If one tries to construct such Riemannian manifolds 
using analogous methods to Hedlund's,
one has to prove a group theoretical property for the fundamental group.
The groups having this property will be called 
\me{groups of bounded minimal generation}.
Any finitely generated abelian group is of 
bounded minimal generation, and every 
group of bounded minimal generation is virtually nilpotent, i.e.\ it 
has a nilpotent subgroup
of finite index. 
Unfortunately, there are only few non abelian groups
that are known to be of bounded minimal generation, e.g. discrete 
subgroups of Heisenberg groups (see section~\ref{secbmg}). 
So this type of construction seems to fail for general fundamental groups.

Therefore we will use a different method that will give us examples for
any finitely generated nilpotent fundamental group.

The Riemannian manifolds we construct 
have a universal covering $\tiM=G\times \ti S$ where
$G$ is a nilpotent Lie-group and $\ti S$ is a 
simply-connected compact manifold.
The commutator group $[G,G]$ acts isometrically on $\tiM$ 
via left multiplication on the first component.
In analogy to the Hedlund metrics on tori we will find two 
types of minimal geodesics on $M$: the \me{left-translated-periodic type} 
and the \me{connection type}.
There are $b_1$ periodic minimal geodesics $c_1,\dots,c_{b_1}$ on $M$
with lifts $\ti{c}_1,\dots,\ti{c}_{b_1}$ on $\tiM$ 
with the following property:
\begin{quote}
  for every
  minimal geodesic $c$ of left-translated-periodic type we can find a lift 
  $\ti{c}$ to $\tiM$ and $i\in\{1,\dots,b_1\}$, $a,b\in\reell$, $g\in[G,G]$ 
  with 
     $$ \ti{c}(t)=g\cdot \ti{c}_i(at+b). $$
\end{quote}
On the other hand every $c\colon\reell\to M$ 
satisfying this property is a minimal geodesic.

A minimal geodesic of connection type will always be a homoclinic or 
heteroclinic connection between two geodesics of left-translated-periodic 
type.

Additionally, the main theorem for this construction 
(Theorem \ref{psymm}) is useful for 
other applications.
For example, we will be able to determine all minimal geodesics 
on nilmanifolds with left-invariant metrics.
Here minimal geodesics are exactly the horizontal lifts
of straight lines on the associated (flat) Jacobi variety $T^{b_1}$.

\head{Acknowledgement}\\
This paper resulted from the author's Diplomarbeit at the University
of Freiburg, Germany. The author wants to thank his supervisor
Victor Bangert for many interesting discussions.
\komment{
and many good ideas 
how to present the results in this paper.
The author profited very much from V.~Bangert's good overview over
existence properties of minimal geodesics and related questions. 
}

\section{$N$-leftinvariant metrics}

In this section we will look at metrics on $G\times \ti S$ as above.
Applying theorem \ref{psymm}, we will be able to reduce the 
classification problem of minimal geodesic on $G\times \ti S$
to the classification of minimal geodesics on $\reell^{b_1} \times \ti S$,
or to be more precise:
\begin{quote} 
the minimal geodesics on $G\times \ti S$ are
exactly the horizontal lifts 
of minimal geodesics  on $\reell^{b_1} \times \ti S$ via the canonical 
Riemannian submersion  
  $$G\times \ti S\to\reell^{b_1} \times \ti S.$$
\end{quote}

We will formulate the theorem in a more general setting.

Let $G$ be a simply connected, nilpotent Lie-group.
The Lie-group exponential map $\exp$ 
is a global diffeomorphism from the Lie-algebra 
$\altg$ 
to the Lie-group $G$, and the Formula of Baker, Campbell and Hausdorff 
states that the pullback of the multiplication on $G$ is a Lie-bracket 
polynomial on $\altg$. Connected subgroups of $G$ correspond 
to Lie-subalgebras of $\altg$ and are therefore 
closed subsets. 
Normal connected subgroups correspond to ideals of $\altg$.
For details and further results about nilpotent Lie-groups 
look for example at \cite{corgreen} and \cite{raghu}.

We fix now a normal connected subgroup $N$ of $G$ 
with Lie-algebra $\altn\subset \altg$.
We will assume that $N$ is contained in the commutator group $[G,G]$ of $G$.
This is equivalent to the condition that $\altn$ is in 
the commutator Lie-algebra $[\altg,\altg]$ of $\altg$.

Let $S$ be a compact manifold;
here we do not assume that $S$ is simply connected but we want $\pi_1(S)$ 
to be finite. 
$G$ acts on $G\times S$ via left multiplication on the first 
component.
Now we take a Riemannian metric $\leerklam_G$ on $G \times S$ 
that is \me{$N$-leftinvariant},
i.e. the subgoup $N$ of $G$ acts isometrically.
Then there is a unique Riemannian metric $\leerklam_{G/N}$ 
on $\left(G/N\right)\times S$ such that 
the canonical projection $p\colon G\times S \to \left(G/N\right)\times S$ 
is a Riemannian submersion.
Vice versa, for every Riemannian metric on
$(G/N)\times S$ there is a (non unique) $N$-leftinvariant 
metric on $G\times S$ such that $p$ is Riemannian. 
 
Additionally we suppose that $\leerklam_G$ is bi-Lipschitz to 
a left invariant metric 
$\leerklam_l$ on $G\times S$, i.e.\ there are constants $c_1,c_2>0$ with
  $$c_1\la v, v\ra_l\leq\la v,v \ra_G\leq c_2\la v,v \ra_l\quad 
    \forall v \in T(G\times S).$$ 
This condition is independent of the choice of the left invariant metric 
$\leerklam_l$.

A vector $v \in T_x(G\times S)$ is called \me{horizontal} if 
$v\perp \kernel T_x p$.

\begin{theorem}\label{psymm}
Let $N$ be a normal connected subgroup of the simply connected, 
nilpotent Lie-group $G$ with $N\subset [G,G]$.
We suppose that $G\times S$ carries a Riemannian metric that
is $N$-left-invariant and bi-Lipschitz to a $G$-left-invariant metric 
and that 
  $$p\colon G\times S \to {G\over N}\times S$$
is a Riemannian submersion.
Then $c\colon \reell \to G\times S$ is a  minimal geodesic on 
$(G\times S,\leerklam_G)$ if and only if
\begin{enumerate}
\item $\dot{c}(t)$ is horizontal for all $t\in \reell$\ \ and
\item $p \circ c\colon \reell \to \left(G/N\right)
      \times S$ is a minimal geodesic on 
      $\left(\left(G/N\right)\times S,\leerklam_{G/N}\right)$.
\end{enumerate}
\end{theorem}

\begin{figure}
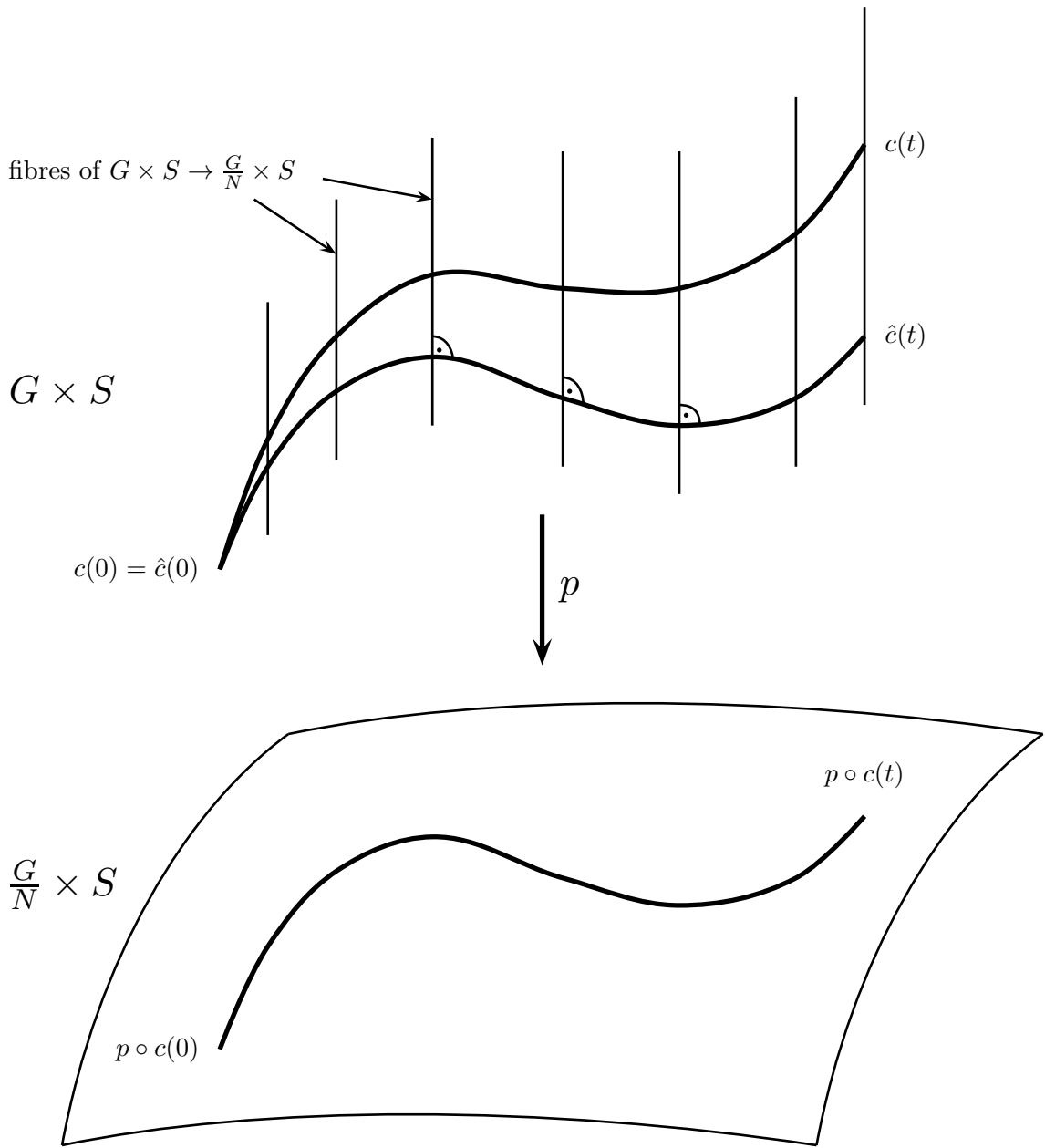

\psset{unit=10mm,showpoints=false}
\pspicture(2,0)(14,18)


\psset{linewidth=2pt}
\def\basecurve{{\pscurve(3.5,3.2)(4.2,4.7)(5.2,5.8)(6.6,6.3)(8.5,5.7)%
(10.2,5.3)(11.9,5.7)(12.9,6.6)}}

\def\basecurveb{{\pscurve(3.5,3.2)(4.2,5.1)(5.2,6.6)(6.6,7.5)(8.5,7.3)%
(10.2,7.3)(11.9,8.1)(12.9,9.4)}}

           \basecurve
           \rput[r](3.2,3.2){$p\circ c(0)$}
           \rput(12.9,7.2){$p \circ c(t)$}
\rput(0,7){\rput[r](3.2,3.2){$c(0)=\hat{c}(0)$}
           \rput[l](13.2,6.6){$\hat{c}(t)$}
           \rput[l](13.2,9.4){$c(t)$}
           \basecurve
           \basecurveb}

\psset{linewidth=1pt}
\rput{0}(0,7)
{
\psline(4.2,3.7)(4.2,7.1)
\psline(5.2,4.8)(5.2,8.6)
\psline(6.6,5.3)(6.6,9.5)
\psline(8.5,4.7)(8.5,9.3)
\psline(10.2,4.3)(10.2,9.3)
\psline(11.9,4.7)(11.9,10.1)
\psline(12.9,5.6)(12.9,11.4)
\psline{->}(5.0,8.9)(6.6,8.6)
\psline{->}(4,8.6)(5.2,7.8)
\rput[r](4.6,9){fibres of $G\times S\to {G\over N}\times S$}
}

\psset{linewidth=\anglestrength}
\rput(0,7)
{
\psbezier(6.6,6.6)(6.8,6.6)(6.9,6.4)(6.9,6.25)
\psdots[linewidth=\angledotstrength](6.7,6.4)
\psbezier(8.5,6.0)(8.7,6.0)(8.8,5.8)(8.8,5.6)
\psdots[linewidth=\angledotstrength](8.6,5.8)
\psbezier(10.2,5.6)(10.4,5.6)(10.5,5.5)(10.5,5.3)
\psdots[linewidth=\angledotstrength](10.3,5.45)
}

\psset{linewidth=2pt}
\rput{0}(8.2,11.0)
{
\psline{->}(0,0)(0,-2.2)
\rput(.4,-1.1){\scalebox{1.5}{$p$}}
}

\rput(1.2,5.6){{\scalebox{1.5}{${G\over N}\times S$}}}
\rput(1.2,12.8){\scalebox{1.5}{$G\times S$}}

\psset{linewidth=1pt}
\psbezier(1.2,1.8)(1.7,4.5)(3.0,6.7)(4.5,7.8)
\psbezier(12.2,1.8)(12.7,4.5)(14.0,6.7)(15.5,7.8)
\psbezier(1.2,1.8)(4.1,2.4)(8.2,2.4)(12.2,1.8)
\psbezier(4.5,7.8)(7.4,8.4)(11.5,8.4)(15.5,7.8)

\endpspicture

\caption{Minimal Geodesics are horizontal}
\end{figure}

\begin{pf}
Let $Z_1(G):=\{x\in G\,|\,xyx^{-1}y^{-1}=e \,\forall y\in G\}$ be 
the center of $G$ and define inductively
$Z_{i+1}(G):=\{x\in G\,|\,xyx^{-1}y^{-1}\in Z_i(G)\,\forall y \in G\}$.

We will prove the theorem for the case $N\subset Z_1(G)$. 
By a straightforward induction on $i$ we then get the theorem for
$N\subset Z_i(G)$ and therefore the general case.

To prove ``$\Leftarrow$'' we suppose that $c$ is not minimal.
If $\dot{c}(t)$ is horizontal for all~$t$,
then $\res{c}{[s,t]}$ has the same length 
as $\res{p\circ c}{[s,t]}$ and therefore
$p\circ c$ cannot be minimal.

For ``$\Rightarrow$'' we suppose 
that $c\colon \reell \to G\times S$ is a minimal geodesic, parametrized by
arclength.
Without loss of generality we can assume that $S$ is simply connected.

For any $v\in \altg$ let $r_{\exp v}$ be the right-translation 
of $G\times S$ by $\exp v$ acting trivially on $S$. 
Then $\vf{v}:=\ddal r_{\exp \al v}$ is a left-invariant vector field 
with vanishing $S$-component.

If $n\in \altn$ then $r_{\exp \al n}$ acts isometrically on $G\times S$,
since $N\subset Z_1(G)$, 
so $\vf{n}$ is a Killing field. 
Noether's theorem (\cite{arnold} 4.20)
implies that 
   $$ \PP_n:=\langle\dot{c}(t),\vf{n}(c(t))\rangle_G $$
is constant in $t$.
We argue by contradiction to show that $\PP_n=0$.

We write 
    $$ \dot{c}(t)=\lambda(t) \vf{n}(c(t)) + c_{\perp}(t) $$
with $c_{\perp}(t)\perp \vf{n}(c(t))$. 

Let $\|\,.\, \|_G$ be the norm of tangential vectors induced by 
$\leerklam_G$. 
If we assume that $\PP_n\neq 0$ we can use the Lipschitz constants between 
$\leerklam_G$ and a left-invariant metric 
to obtain
constants $\Konst_1,\Konst_2>0$ in the inequalities:

  $$ \left|\lambda(t)\right|={\left|\Big\langle\dot{c}(t),
      \vf{n}(c(t))\Big\rangle_G \right|\over
         \left\| \vf{n}(c(t))\right\|_G^2}
      = {|\PP_n|\over \left\| \vf{n}(c(t))
         \right\|_G^2}
         < \Konst_1$$
  $$\|\lambda(t)  \vf{n}(c(t))\|_G 
\komment{
      & = & {\left|\Big\langle\dot{c}(t),
      \vf{n}(c(t))\Big\rangle_G \right|\over
      \left\| \vf{n}(c(t)) \right\|_G} \\
}
       =  {|\PP_n|\over \left\| \vf{n}((c(t))
        \right\|_G} 
       >  \Konst_2>0$$

Then the curve $\hat{c}$ defined by
  $$\hat{c}(t):= c(t) \cdot \exp \bigg(\Big(\int_0^t 
                       - \lambda(t')dt'\Big)n\bigg)$$
satisfies
$p \circ \hat{c} = p\circ c$,
$\hat{c}(0)=c(0)$
and 
$\dot{\hat{c}}(t) \perp \vf{n}$ (see also Figure 1).
After identification of $T_{\hat{c}(t)}M$ and $T_{c(t)}M$ via left translation,
$\dot{\hat{c}}(t)$ is equal to $c_{\perp}(t)$.

So we know that $\|\dot{\hat{c}}\|_G\leq\sqrt{1-{\Konst_2}^2}$.
Writing $d_G$ for the distance induced by $\leerklam_G$ we obtain
\begin{equation}\label{wider}
  t=d_G(c(0),c(t))  \leq 
  \underbrace{d_G(\hat{c}(0),\hat{c}(t))}_{\leq 
   \left(\sqrt{1-{\Konst_2}^2}\right) t}
  \: + \: d_G(\hat{c}(t),c(t)) \quad \forall t>0.
\end{equation}
We use a result of Pansu (\cite{pansu}) to state that there is a 
constant $\Konst_3(n)$ not depending 
on $\al>0$, $s\in S$ and $g \in G$ such that
\begin{equation}\label{quadwachs}
  d_G\big((g \exp \al n,s),( g,s)\big) \leq \Konst_3(n)\,(\sqrt{\al}+1). 
\end{equation}

Pansu did not prove exactly this statement, but
the proof of it is completely analogous to the proof of  
\cite{pansu} no.\ (62) if we use the fact that $\exp \al n$
is in the commutator group. Another proof 
using more elementary methods can be found in
\cite{ammann}.
  
Together with $|\lambda (t)|<\Konst_1$ 
inequality (\ref{quadwachs}) 
contradicts (\ref{wider}), so we get $\PP_n=0$ for every $n$ in the 
Lie algebra of $N$, i.e.\ $\dot{c}(t)$ is horizontal.
This implies that $p\circ c$ is parametrized by arclength. 

\begin{figure}
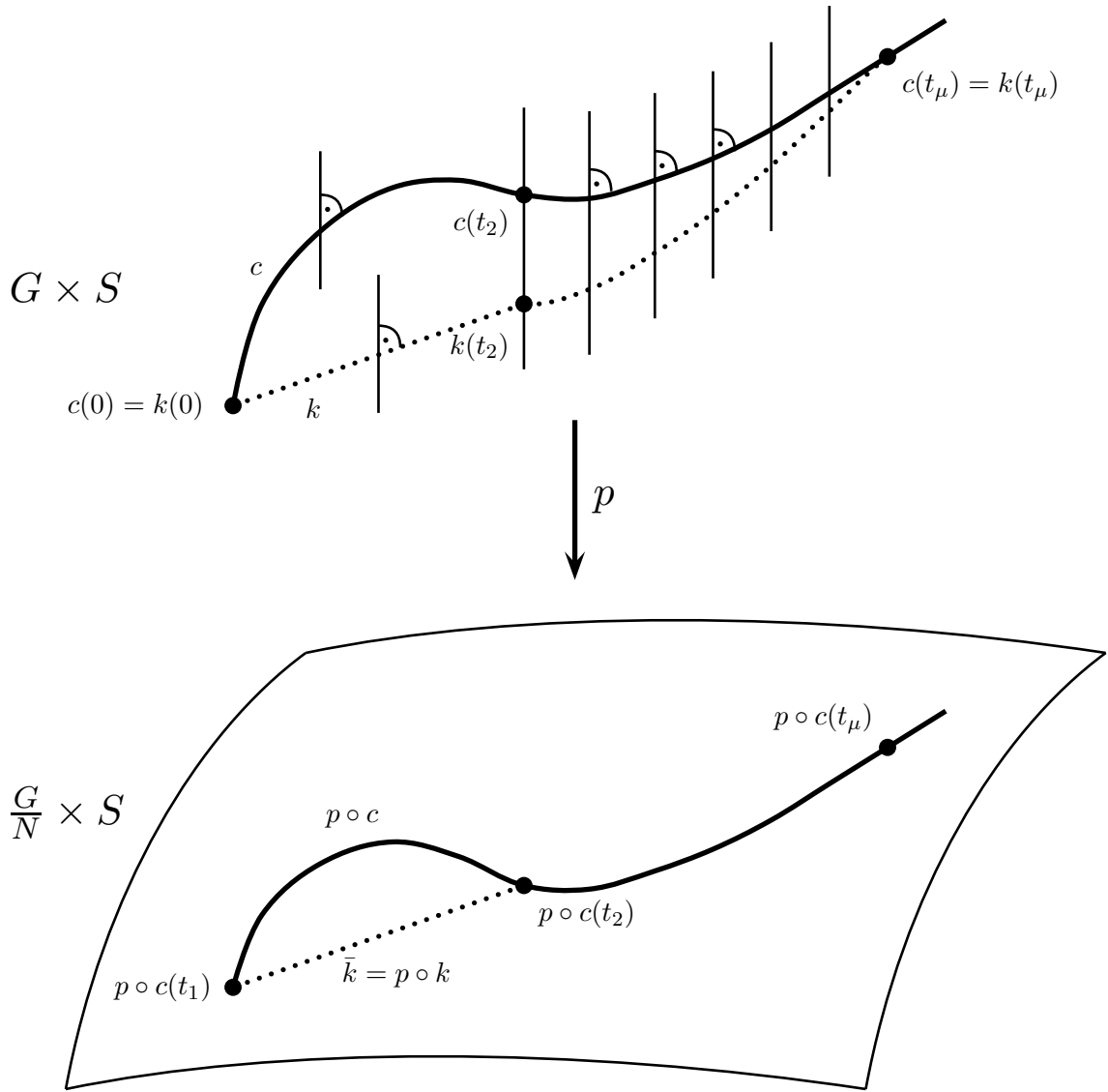

\psset{unit=10mm,showpoints=false}
\pspicture(2,0)(14,18)

\psset{linewidth=\curvestrength}
\def\basecurve{{\pscurve(3.5,3.2)(3.9,4.2)(4.7,4.9)(5.7,5.2)(6.6,5.0)(7.5,4.6)%
(8.4,4.55)(9.3,4.8)(10.1,5.1)(10.9,5.5)(11.7,6.0)(12.5,6.5)(13.3,7.0)}}

\def\basecurvea{{\pscurve(3.5,3.2)(3.9,4.6)(4.7,5.6)(5.7,6.2)(6.6,6.3)(7.5,6.1)%
(8.4,6.05)(9.3,6.3)(10.1,6.6)(10.9,7.0)(11.7,7.5)(12.5,8.0)(13.3,8.5)}}

\def\basecurveb{{%
  \psset{linewidth=\curvestrength,linestyle=dotted}
  \psline(3.5,3.2)(7.5,4.6)
  \psecurve(6.6,5.0)(7.5,4.6)(8.4,4.8)(9.3,5.3)(10.1,5.85)%
           (10.9,6.5)(11.7,7.25)(12.5,8)(13.3,8.5)
  }}

       \basecurve          
           \rput[rb](12.3,6.7){$p \circ c(t_{\mu})$}
           \rput[r](3.2,3.2){$p\circ c(t_1)$}
           \rput[lt](7.7,4.4){$p\circ c(t_2)$}
           \rput[rb](5.5,5.4){$p \circ c$}
           \psdots[linewidth=\dotstrength](3.5,3.2)
           \psdots[linewidth=\dotstrength](7.5,4.6)
           \psdots[linewidth=\dotstrength](12.5,6.5)
           \psline[linestyle=dotted,linewidth=\curvestrength](3.5,3.2)(7.5,4.6)
           \rput[lt](5,3.6){$\bar{k}=p \circ k$}
\rput(0,8){\rput[r](3.1,3.2){$c(0)=k(0)$}
           \rput[lt](12.7,7.8){$c(t_{\mu})=k(t_\mu)$}
           \psdots[linewidth=\dotstrength](3.5,3.2)
           \psdots[linewidth=\dotstrength](7.5,6.1)
           \psdots[linewidth=\dotstrength](12.5,8.0)
           \psdots[linewidth=\dotstrength](7.5,4.6)
           \rput[rt](7.3,5.9){$c(t_2)$}
           \rput[rt](7.3,4.2){$k(t_2)$}
           \rput[rb](3.9,5.0){$c$}
           \rput[lt](4.5,3.3){$k$}
           \basecurvea
           \basecurveb
}

\psset{linewidth=\fibrestrength}
\rput{0}(0,8)
{
\psline(7.5,3.7)(7.5,7.3)
\psline(8.4,3.9)(8.4,7.25)
\psline(9.3,4.4)(9.3,7.5)
\psline(10.1,4.95)(10.1,7.8)
\psline(10.9,5.6)(10.9,8.2)
\psline(11.7,6.35)(11.7,8.7)
\psline(4.7,4.8)(4.7,6.7)
\psline(5.5,3.1)(5.5,5.0)
}

\psset{linewidth=\anglestrength}
\rput{0}(0,8)
{
\psbezier(8.4,6.45)(8.6,6.45)(8.7,6.35)(8.7,6.15)
\psdots[linewidth=\angledotstrength](8.5,6.25)
\psbezier(9.3,6.7)(9.5,6.7)(9.6,6.6)(9.6,6.4)
\psdots[linewidth=\angledotstrength](9.4,6.5)
\psbezier(10.1,7.0)(10.3,7.0)(10.4,6.9)(10.4,6.7)
\psdots[linewidth=\angledotstrength](10.2,6.8)
\psbezier(4.7,6.1)(4.9,6.1)(5.0,6.0)(5.0,5.8)
\psdots[linewidth=\angledotstrength](4.8,5.9)
\psbezier(5.5,4.3)(5.7,4.3)(5.8,4.2)(5.8,4.0)
\psdots[linewidth=\angledotstrength](5.6,4.1)
}

\psset{linewidth=\curvestrength}
\rput{0}(8.2,11.0)
{
\psline{->}(0,0)(0,-2.2)
\rput(.4,-1.1){\scalebox{1.5}{$p$}}
}

\rput(1.2,5.6){{\scalebox{1.5}{${G\over N}\times S$}}}
\rput(1.2,12.8){\scalebox{1.5}{$G\times S$}}

\psset{linewidth=\framestrength}
\psbezier(1.2,1.8)(1.7,4.5)(3.0,6.7)(4.5,7.8)
\psbezier(12.2,1.8)(12.7,4.5)(14.0,6.7)(15.5,7.8)
\psbezier(1.2,1.8)(4.1,2.4)(8.2,2.4)(12.2,1.8)
\psbezier(4.5,7.8)(7.4,8.4)(11.5,8.4)(15.5,7.8)

\endpspicture

\caption{Minimal Geodesics project to Minimal Geodesics}
\end{figure}

It remains to show that $p\circ c$ is minimal. In order to
prove it we assume the opposite, i.e.\  
  $$\De:=t_2-t_1- d_{G/N}(p\circ c (t_1),p\circ c(t_2))>0.$$
Now take a shortest geodesic $\bar{k}\colon[t_1,t_2]\to {(G/N)}\times S$ 
from $p \circ c(t_1)$ to $p \circ c(t_2)$ (see also Figure 2). 
This shortest geodesic has a unique horizontal lift
$k\colon [t_1,t_2]\to G\times S$ with $c(t_1)=k(t_1)$.

As the Lie exponential map is a diffeomorphism
there is a unique $n\in \altn$, such that
$k(t_2)=c(t_2)\cdot \exp n$.
For $\mu>0$ we now extend $k$ continuously by
  $$k(t):=c(t)\cdot \exp([1-\mu(t-t_2)]n) 
    \qquad t_2 \leq t \leq t_{\mu}:=t_2 + 1/\mu .$$
So $k\colon{[t_1,t_{\mu}]}\to G\times S$ is also a curve from $c(t_1)$ to
$c(t_{\mu})$.
We will prove that $k$ is shorter 
than $\res{c}{[t_1,t_{\mu}]}$ for small $\mu>0$.
There is a unique $c':\reell \to \altn^{\perp} \subset \altg$ 
such that the $G$-component
of $\dot{c}(t)$ is $\vf{c'(t)}(c(t))$. 
On $(t_2,t_{\mu})$ the $G$-component of $\dot{k}(t)$ is
  $$\vf{c'(t)}(k(t))- \mu \vf{n}(k(t)),$$
whereas the $S$-components of $\dot{k}(t)$ and $\dot{c}(t)$ are equal up to 
left (or right) translation.

So as $c$ is horizontal
\begin{eqnarray*}
 \|\dot{k}(t)\|_G &   =   & \sqrt{
                        \|\dot{c}(t)\|_G^2
                        + \mu^2 \| \vf{n}(k(t)) \|_G^2}\cr
                  & \leq  & 1+ {1\over 2}\mu^ 2 \| \vf{n}(k(t))\|_G^2. 
\end{eqnarray*}
For $\mu$ small enough we get
  $$\LL\big(\res{c}{[t_1,t_{\mu}]}\big) -\LL(k)  \geq
            \Delta - {1\over \mu}
    {\mu^2\over 2}\sup_{g\in G\atop s\in S}\|\vf{n}(g,s)\|_G^ 2>0,$$
which contradicts the minimality of $c$. Therefore $p\circ c$
is minimal.\qed
\end{pf}

Using Theorem \ref{psymm} we now know the minimal geodesics
on any nilpotent Lie-group with a left-invariant metric.

\begin{corollary}
Let $G$ be a nilpotent Lie-group
with a left-invariant metric.
The minimal geodesics on
$G$ are exactly the curves of the form 
  $$c(t)=g\cdot \exp tv$$
with $g\in G$ and $v\in \altg$ and $v\perp [\altg,\altg]$.
\end{corollary}

\remark A similar type of orthogonality relation was discovered by 
Patrick Eberlein, Ruth Gornet and Dorothee Sch\"uth when they investigated
the following problem. A geodesic $c$ on a nilpotent Lie-group $G$ with 
left-invariant metric is called \me{periodic} if there are
$g\in G$, $\lambda>0$ with $c(t+\lambda)=g\cdot c(t)\,\,\forall t$.
A~necessary condition for periodicity is that periodic geodesics are 
orthogonal to certain terms built by commutators 
(\cite{eberlein}Cor.~4.4, \cite{gornet}3.1). 
\bigskip

Now we try to lift asymptotic behavior from $(G/N)\times S$ to
$G\times S$. Here we have to pay attention to the following fact:
if $\dim N>0$ then there are curves  
$\bar{\ga}_1,\bar{\ga}_2\colon \reell \to (G / N)\times S$
that are asymptotic to each other but do not have 
horizontal lifts that are asymptotic to each other.

The situation is different if we replace ``asymptotic'' 
by ``exponentially asymptotic''.

\begin{definition}\label{expasymptdef}
The (parametrized) curves $\ga_1$ and $\ga_2$ 
are \me{exponentially asymptotic}
for $t\to \infty$ 
if there are constants $\Konst_a, \Konst_b >0$ such that 
  $$d\left(\ga_1(t),\ga_2(t)\right)<\Konst_a e^{-\Konst_b t} 
    \mbox{ for large } t .$$
\end{definition}
Here $d$ is the Riemannian distance. This definition is
invariant under bi-Lipschitz change of the metric.
The definition of exponentially asymptotic for $t\to -\infty$ is analogous.

\begin{proposition}\label{expasymptlift}
Let $G\times S$ carry a Riemannian metric that is $N$-left-invariant and 
invariant under left-action of a lattice $\Gamma$ of $G$. Choose a Riemannian
metric on the quotient such that 
$p:G\times S \to (G/N)\times S$ 
becomes a Riemannian submersion.
Furthermore let $\ga_1,\ga_2\colon \reell \to G\times S$ 
be piecewise $C^1$-curves with  
bounded $\|\dot{\ga}_i\|_G$ and horizontal with respect to $p$. 
Then $p\circ \ga_1$ is exponentially asymptotic to $p \circ \ga_2$ if 
and only if there is an $n_{\infty}\in N$ such that $n_{\infty}\cdot\ga_1$
is exponentially asymptotic to $\ga_2$.
\end{proposition}

\spacebetweentheoremandproof

\begin{pf} 
We only have to prove the ``only if''. And it is sufficient 
to prove this for $N \subset Z_1(G)$ as the 
general case follows by induction.
The cases $t \to + \infty $ and $t \to -\infty$ are totally symmetric,
so the case $t \to -\infty$ will be omitted.

Because of the action of $\Gamma$ and the compactness of $S$ 
the injectivity radius of ${G/N} \times S$ is positive.

At first we can choose $t_0$ with 
  $$d(p \circ \ga_1(t), p \circ \ga_2(t)) < {1\over 2} 
                    \injrad ({G\over N}\times S) \qquad \forall t \geq t_0.$$
We glue a surface $A\colon [1,2]\times [t_0,\infty)\to (G/N)\times S$ 
between the $\res{p\circ \ga_i}{[t_0,\infty]}$ 
such that $A(i,.)=p \circ \ga_i\quad (i=1,2)$ and $A(.,t)$ is the shortest 
curve from $p\circ \ga_1(t)$ to $p \circ \ga_2(t)$.

As $p \circ \ga_1(t)$ and $p \circ \ga_2(t)$ are exponentially asymptotic,
there are constants $\Konst_a, \Konst_b>0$ with 
  $$\area \res{A}{[1,2]\times [t_1,t_2]} < \Konst_a e^{-\Konst_b t_1} \quad 
    (t_0\leq t_1\leq t_2 \leq \infty).$$

Now lift $A$ to $\ti{A}\colon [1,2]\times [t_0,\infty) \to G \times S$ such
that $\ti{A}(1,.)=\ga_1$ and $\ti{A}(.,t)$ is horizontal $\forall t$.
There is an $n\colon [t_0,\infty) \to N$ with 
  $$\ti{A}(2,t) \cdot n(t) = \ga_2(t).$$
In general $n$ is non-constant and therefore $\ti{A}(2,.)$ is non-horizontal,
but we will show that $n(t)$ converges for $t\to \infty$.

Note that $G\times S\to (G/N)\times S$ is a principal $N$-bundle.
As $N$ acts isometrically, the horizontal planes determine 
a connection-1-form $\om \colon T(G\times S) \to \altn$.
(For details on connection-1-forms see \cite{kobayashi}, Chapter II.)
Then $d\om$ is the curvature of the connection and $\|d\om\|$ is uniformly 
bounded on $G\times S$.

As $\ga_i$ and $\ti{A}(.,t)$ are horizontal 
  $$\int_{\res{\ga_i}{[t_1,t_2]}}\om =0 \quad \mbox{and} \quad 
    \int_{\ti{A}(.,t)} \om =0$$
  $$\int_{\partial \left( \res{\ti{A}}{[1,2]\times[t_1,t_2]}\right)}\om
    =  \int_{\res{\ti{A}(2,.)}{[t_1,t_2]}} \om = n(t_1) n(t_2)^{-1}$$
Using Stoke's Theorem we get
\begin{eqnarray*}
    d\left(n(t_1)n(t_2)^{-1},e\right)\leq \int_{\res{\ti{A}}{[1,2]\times [t_1,t_2]}}\|d\om\|
    & \leq  & \area \left( \res{A}{[1,2]\times [t_1,t_2]}\right) 
              \sup_{G\times S} \|d\om\|\cr
    & \leq  & \sup \|d\om\| \Konst_a e^{-\Konst_b t_1}.
\end{eqnarray*}
So $n(t)$ converges exponentially 
and $n_{\infty}:=\lim_{t\to \infty} n(t)$ gives the proposition.
\qed
\end{pf}

\section{Hedlund examples}

In this chapter we will give a slight generalisation of the Hedlund 
examples presented by Bangert in \cite{bangertmin}, section 5.
The proofs are only small 
variations of Bangert's proofs, so we
will skip them.

In this section we construct similar metrics, which we will also call
``Hedlund metrics''. These metrics are defined on manifolds of the form 
$M=T^{b_1}\times S$ where 
$T^{b_1}=\reell^{b_1}/\Z^{b_1}, b_1\geq 1$ is the torus and $S$ is 
an arbitrary compact connected
manifold with  finite $\pi_1(S)$.
We exclude the case $M=T^2$ by assuming  
$\dim S>0$ or $b_1\neq 2$.
Note that $b_1$ is the first Betti number of $M$.

We denote the standard flat metric on $T^{b_1}$ by $\leerklam_T$
and we choose a metric $\leerklam_S$ on $S$. 
The product metric on 
$M$ will be called $\leerklam_{T\times S}$.
The vectors of the canonical basis $e_1,\dots, e_{b_1}$ of $\reell^{b_1}$
induce $\leerklam_{T\times S}$-orthonormal 
vector fields $E_1,\dots,E_{b_1}\colon M \to TM$. 
 
The Hedlund metrics $\leerklam_H$
will be defined in Definition \ref{hedldef}. 
In this definition we use $b_1$ closed curves 
$c_1,\dots,c_{b_1}$ 
on $M$ that will become the only geodesics that are minimal and closed.
To define them, we have to distinguish two cases.

In the case ``$b_1\neq 2$''
choose $s\in S$ and define 
  $$c_i(t)\,:=\, \left({{1\over 2}e_i + t e_{i+1}\over \Z^{b_1}},\,s\,\,\right) $$
for all $t\in \reell$, where $e_{b_1+1}:=e_1$.
 
In the case ``$b_1=2$'' we have assumed that $\dim S >0$, so we can choose 
different $s_1,s_2\in S$ and define
  $$c_i(t)\,:=\, \left({te_i\over \Z^{b_1}},\,s_i\,\right)$$
for all $t\in \reell$.

In both cases let $L_i$ be the trace of $c_i$. 
The fact that $L_i$ and $L_j$ are disjoint for $i\neq j$ plays 
an important role in the proofs. 
The construction of Hedlund type metrics on $T^2$ fails because
such $c_i$ and $L_i$ do not exist on $T^2$.    

Now define 
$U_{\ep}(L_i)$ to be the $\ep$-neighborhood of $L_i$
with respect to $\leerklam_{T\times S}$.
 
For $\ep>0$ (that will be chosen very small) we define
in analogy to Definition 5.1 of \cite{bangertmin}:
\begin{definition}                    \label{hedldef}
 $\leerklam_H$ is an $\ep$-Hedlund metric on $M$
iff there are $\ep_1,\dots,\ep_{b_1}\in (0,\ep]$ such that
for $i=1,\dots, b_1$:

$\matrix{
(P1) & \la v,v \ra_H \leq (1+\ep)^2\la v,v \ra_{T\times S} \hfill 
        & \forall v \in TM\hfill\cr
(P2) & \la E_i(x),E_i(x)\ra_H = \ep_i^2 \hfill         
        & \forall x \in L_i\hfill\cr
     & \la v,v \ra_H \geq \ep_i^2\la v,v\ra_{T \times S} \hfill  
        & \forall v\in T_xM, x\in  L_i \hfill\cr
     & \la v,v \ra_H > \ep_i^2\la v,v\ra_{T \times S} \hfill 
        & \forall v\in T_xM\ohne\{0\}, x\in U_\ep(L_i)\ohne L_i\hfill\cr 
(P3) & \la v,v \ra_H \geq  \la v,v\ra_{T \times S} \hfill
        & \forall v\in T_xM, x\not\in \bigcup_j U_\ep(L_j) \hfill\cr
}$
\end{definition}
The following propositions \ref{beschrlaeng},
\ref{beschrtub} and \ref{geodasympt} are analogues to 
Propositon 5.2, Proposition 5.3 and Corollary 5.4 of \cite{bangertmin}.
Because of the definition of $\ep$-Hedlund metric it is clear
that any statement in these propositions that holds for $\ep>0$ also holds
for any $\ep'\in(0,\ep)$. 

\begin{proposition}\label{beschrlaeng}
There is an $\ep>0$ and a $\Konst_1\in \reell$ such that
for any $\ep$-Hedlund metric on $M$ and any arclength-parametrized 
minimal geodesic $c$ with respect 
to this metric the length 
of
  $$A:=c^{-1}\big(M\ohne \bigcup_i U_{\ep}(L_i)\big)\subset \reell $$
is bounded by $\Konst_1$.
\end{proposition}
That means that $c$ ``stays out of $\bigcup_i U_{\ep}(L_i)$ only 
for a bounded time''.
As an immediate consequence $c$ cannot change its ``tube'' too often.
To make this precise we define:

\begin{definition}\label{changenumber}
Let $\ep>0$ be so small that the $U_{\ep}(L_i)$ $(i=1,\dots , b_1)$
are disjoint and let $c\colon \reell \to M$ be a minimal geodesic.
We define the change number $\CC(c)\in \natur \cup\{0,\infty\}$ 
to be the supremum of all $n \in \natur$ such that we find  
$t_0<t_1<\dots<t_n$ and $i_j \in \{1,\dots,b_1\}$
with $c(t_j)\in U_{\ep}(L_{i_j})$ and $i_j\neq i_{j+1}$.
\end{definition}

\begin{proposition}\label{beschrtub}
There is an $\ep>0$ and $\Konst_2\in \natur$ such that 
$\CC(c)\leq\Konst_2$ for any minimal geodesic
$c$ with respect to any $\ep$-Hedlund metric on $M$.
\end{proposition}

\remark
If $S=\{\mbox{one point}\}$ Bangert proved in \cite{bangertmin}
that we can even find 
$\ep>0$ with $\Konst_2:=b_1$. For general $S$ this statement does not hold.
\smallskip

\begin{proposition}\label{geodasympt}
There is an $\ep>0$ such that every minimal geodesic on an 
$\ep$-Hedlund metric 
on $M$ is asymptotic in each of its senses 
to one of the $L_i$'s.
\end{proposition}

In section \ref{mainconst} we will need a stronger version, 
so we formulate a supplement.

\begin{supplement}\label{zusatz}
If the $c_i$ are even hyperbolic closed geodesics, e.g.\ if
  $$A_{jk}(x):=E_k\big(E_j\big(\la E_i,E_i\ra_H\big)\big)(x) 
    \qquad{j,k\neq i}$$
is positive definite for all $x\in L_i$, 
then any minimal geodesic $c$ 
is exponentially asymptotic to one of the $L_i$ in each of
its senses (see Definition \ref{expasymptdef}).
\end{supplement}

\remark
It is also possible to formulate analogues to the propositions 
5.6, 5.7 and 5.8 from \cite{bangertmin}.
\smallskip

\section{Lattices in nilpotent Lie-groups}

Here we will summarize some facts used in the next
section. For the discrete group or Lie-group $G$ we define 
the descending central series $(G^i)_{i\in \natur}$
inductively by $G^1:=G$ and $G^{i+1}:=[G,G^i]$.
Then $G$ is nilpotent iff $G^i=\{e\}$ for sufficiently big 
$i\in \natur $.

\begin{theorem}[Malcev, \cite{raghu} theorem 2.18]\label{thmalcev}
A group $\Ga$ is isomorphic to a lattice in 
a nilpotent, simply connected Lie-group iff $\Ga$ is finitely generated,
nilpotent and torsion free.
\end{theorem}

\begin{theorem}\label{faktorgitter}
Let $\Ga$ be a lattice in the nilpotent, simply connected 
Lie-group $G$ and $N$ a closed normal subgroup (not necessarily connected),
$p\colon G \to G/N$.\newline
If two of the following three conditions are true, the third follows:
\begin{enumerate}
\item $\Ga\cap N$ is a lattice in $N$,
\item $p(\Ga)$ is a lattice in $G/N$,
\item $\Ga$ is a lattice in $G$.
\end{enumerate}
\end{theorem}
``1.\ and 2.\ $\Rightarrow$ 3.''
and ``1.\ and 3.\ $\Rightarrow$ 2.'' 
are proved in \cite{corgreen} lemma 5.1.4, the proof 
of ``2.\ and 3.\  $\Rightarrow$ 1.'' is straightforward.

\begin{theorem}\label{abstgitter}
Let $\Ga$ be a lattice in $G$, then $\Ga^i$ is a lattice in $G^i$
for $i \in \natur$.
\end{theorem}
This follows from the theory of Malcev bases (\cite{corgreen}) 
and \cite{corgreen} corollary~5.4.5. 
It is a slight generalisation of \cite{raghu} corollary~1 of theorem~2.3
saying that $\Ga \cap G^i$ is cocompact in $G^i$. 

Using the theory of Malcev bases it is also evident that 
  $$\rank {\Ga \over[\Ga,\Ga]}=\dim {G \over [G,G]}.$$

\section{Main construction}\label{mainconst}
In this section we will construct our examples with minimal geodesics in only
``few directions'' by combining 
the results we obtained in the previous sections.

\pagebreak[2]
\begin{theorem}\label{maintheo}
For any finitely generated nilpotent group $\Pi_1$ we find a
connected compact Riemannian manifold $(M,\leerklam_M)$ satisfying:
\begin{description}
\item[\rm (1)] $\pi_1(M) = \Pi_1$
\item[\rm (2)] $M$ has a universal 
      covering $\tiM= G \times \ti{S}$
      where $G$ is a  nilpotent Lie-group and $\ti{S}$ is a
      compact manifold.
\item[\rm (3)] The commutator group $[G,G]$ acts isometrically on the 
      Riemannian covering $\tiM$
      via left multiplication on the first component.
\item[\rm (4)] There are minimal geodesics 
     $c_i\colon \reell\to M$ $(i\in \{1,\dots,b_1\})$ with lifts 
     $\ti{c}_i\colon \reell \to \tiM$ such that 
     every minimal geodesic $\ga\colon \reell \to M$ is of one of the 
     following types:
  \begin{description}
       \item[Type I: left-translated-periodic]\ \\
                      $\ga$ has a lift $\ti{\ga}\colon \reell \to \tiM$ 
                      such that there are $a,b\in \reell$, $g \in [G,G]$,
                      $i \in \{ 1,\dots,b_ 1\}$
                      with $\ti{\ga}(t)=g \cdot c_i (at+b)$.
       \item[Type II: connection type]\ \\
                      $\ga$ is not of Type I, but there are minimal 
                      geodesics $\ga_{+}$ and $\ga_{-}$ of Type I such that 
                      $\ga(t)$ is exponentially asymptotic to 
                      $\ga_{\pm}(t)$ for $t\to \pm \infty$.
   \end{description}
\end{description}
\end{theorem}

Any compact manifold is a finite CW-complex and therefore the 
fundamental group of any compact manifold is finitely genrated.

The author has presented another version of this construction
in \cite{ammann} and proves that if $\Pi_1$ has no torsion and $\Pi_1\neq \Z^2$
we even get an example as above with $\tiM=G\times \{\mbox{one point}\}$
(but without property (3) in the case $b_1=2$). 
It uses the fact that subgroups of the 3-dimensional Heisenberg group
are of bounded minimal generation (Theorem \ref{allebmg}).  
So $T^2$ is the only 
nilmanifold that does not admit a metric of Hedlund type.
  
To prove theorem \ref{maintheo} we will use:

\begin{theorem}[K.A. Hirsch, \cite{baumsl} theorem 2.1]\label{hirsch}
Let $\Pi_1$ be a finitely generated nilpotent group. 
Then $\Pi_1$ can be embedded 
as a subgroup into the direct product $D=\Ga\times F$ where 
$\Ga$ is a torsionfree nilpotent finitely generated group
and $F$ is a finite nilpotent group. 
\end{theorem}

The elements in $\Pi_1$ of finite order form a normal subgroup $T=F\cap \Pi_1$.
Looking at Baumslag's proof of the theorem of Hirsch, we immediately see 
that $\Ga$ can be chosen as $\Pi_1/T$ and that 
$$
\renewcommand{\arraystretch}{5.5}
\hfill
\begin{array}{r@{\hskip 2cm}l}
\rnode{a}{\Pi_1} & \rnode{b}{\Gamma_{}} \times F \\
                 & \rnode{c}{\Gamma_{}}={\Pi_1/T} \\
\end{array}
\psset{nodesep=3pt}
\ncline{>->}{a}{b}
\ncline{->>}{a}{c}
\ncline{->>}{b}{c}
\hfill
$$
commutes.
\bigskip

\begin{pf*}{Proof of Theorem \ref{maintheo}}
\par\noindent
Suppose $\Pi_1$ to be embedded in $D=\Ga \times F$ as above.
 
Now embed $\Ga$ as a lattice in the nilpotent Lie-group $G$ 
(Theorem \ref{thmalcev}).

Abelianisation of the above diagram and tensoring with $\reell$
shows that 
  $$\rank {\Ga \over [\Ga,\Ga]}
    =\rank {\Pi_1\over[\Pi_1,\Pi_1]} \big(=b_1\big),$$ 
so there is a natural isomorphism
  $${G\over [G,G]}\,\ti{\to} \,\reell^{b_1}$$
that maps the image of $\Ga$ to $\Z^{b_1}$.

It is well-known that for $n\geq 4$ 
any finitely presented group
is the fundamental group of an $n$-dimensional compact manifold 
(\cite{massey2} p.~114). 
So let $S$ be a $4$-dimensional compact Riemannian manifold
with $\pi_1(S)=F$.

Now take an $\ep$-Hedlund metric on $T^{b_1}\times S$ with hyperbolic $c_i$ 
and choose a metric on 
$(\Ga\backslash G)\times S$ such that 
  $$(\Ga\backslash G) \times S \to T^{b_1}\times S$$
is a Riemannian submersion and such that $[G,G]$ acts isometrically on 
$G\times S$.
As $(\Ga\backslash G)\times S$ has fundamental group $D$,
we can find a Riemannian covering $M$ with $\pi_1(M)=\Pi_1$.
So (1), (2) and (3) are fulfilled.

Using Theorem \ref{psymm}, Proposition \ref{geodasympt},
Supplement \ref{zusatz} and Proposition \ref{expasymptlift}
we also get (4).
\qed
\end{pf*}
\bigskip

In the remaining part of this section we will discuss some related topics and 
some additional properties of the above examples.

For arbitrary compact manifolds $M$ with $\Pi_1:=\pi_1(M)$ we get 
via the Hurewicz map 
  $$H_1(M,\Z)={\Pi_1\over [\Pi_1,\Pi_1]} \quad \mbox{and} \quad 
    H_1(M,\reell)=H_1(M,\Z)\otimes_{\Z}\reell.$$
Dividing $H_1(M,\reell)$ by the image of $H_1(M,\Z)$ we get a 
torus $T^{b_1}$ in analogy to the above constructions.
This torus is known as the Jacobi variety (\cite{gromov2}4.21).
$H_1(M,\reell)$ also carries a norm, the ``stable norm $\|\,.\,\|$'',
induced by the Riemannian structure of $M$ 
(\cite{federer},\cite{gromov2}4.18,\cite{bangerticm}).
$\left(\tiM,\ep\leerklam\right)$ converges in the Gromov-Hausdorff-sense
to $(H_1(M,\reell),\linebreak[2] \|\,.\,\|)$ 
for $\ep\to 0$ if $\Pi_1$ is abelian.
For nilpotent $\Pi_1$ it converges to a Carnot-Caratheodory space 
(\cite{pansu}) and the stable norm is essential for 
measuring distances on this space.

Bangert used the stable norm to prove an existence theorem for minimal
geodesics (\cite{bangertmin},\cite{bangerticm}). 
As a corollary he proved the existence of at least $b_1$ different geodesics
such that the ``rotation set'' $\RR(c_i)$ of each $c_i$ contains only one 
vector in $H_1(M,\reell)$ and $\bigcup_i \RR(c_i)$ is a basis of 
$H_1(M,\reell)$.
In the above examples the geodesics $c_1,\dots,c_{b_1}$ have the properties
of the geodesics whose existence has been shown by Bangert:
each $\RR(c_i)$ contains only one vector and their union is a basis.
For our examples the stable norm written in this basis is just
\begin{equation}\label{stanor}
  \|x\|=\sum_{i=1}^{b_1} \ep_i |x^i|.
\end{equation}    
Equation (\ref{stanor}) can be seen from Bangert's theorem 
and the characterisation of the minimal geodesics
or just using the fact that if a Riemannian submersion 
$p\colon M_1\to M_2$ of compact manifolds $M_i$ induces 
an isomorphism $p_{\#}\colon H_1(M_1,\reell)\to H_1(M_2,\reell)$,
then $p_{\#}$
preserves the stable norm.

This last fact can also be used to construct metrics on nilmanifolds 
with non-left-invariant metric on the universal covering that
have a smooth unit ball of the stable norm.
(Just take a suitable lift of a flat metric on the Jacobi variety.)

\section{Expressway metrics}

In the previous section we proved that for every 
finitely generated nilpotent group $\Pi_1$
there is a Riemannian manifold $M$ with $\pi_1(M)=\Pi_1$ and only
few directions of minimal geodesics. On the other hand, many properties
concerning minimal geodesics only depend on the fundamental group and
an induced distance on it. Therefore it seems likely that we could find
a suitable Riemannian metric on every compact manifold with nilpotent
fundamental group. 
Moreover the metrics in the last section admit continuous 
families of minimal geodesics if $\dim [G,G]>0$, 
whereas Hedlund's original examples only
admit very few ones. 
So it would be interesting to generalize Hedlund's methods directly.
 
In this section we try to use Hedlund's method 
(generalized by Bangert \cite{bangertmin}) directly to construct
Hedlund type metrics on arbitrary (compact) manifolds with nilpotent 
fundamental group $\Pi_1$. It turns out that we succeed only if $\Pi_1$ has 
an algebraic property that we call ``bounded minimal generation''. We can
prove that lattices in Heisenberg groups have this property but we do not 
know if this is true for all finitely generated nilpotent groups or not.

The following construction seems to be very special, but the author 
thinks that it will be difficult to find a different construction
without getting a problem similiar to the bounded-minimal-generation-problem
described in the next section.
We will omit an explicit definition of the metrics and proofs of 
the statements concerning these metrics as the exact formulae 
only give little insight in what happens. (For details see \cite{ammann}).
Instead we will give an informal description.

\begin{figure}[ht]
\unitlength=0.30mm
\begin{picture}(400.00,230.00)(-8,55)
\fatlines
\bezier{50}(50.00,100.00)(50.00,92.00)(60.00,85.00)
\bezier{100}(60.00,85.00)(76.00,75.00)(100.00,75.00)
\bezier{100}(100.00,75.00)(125.00,75.00)(140.00,85.00)
\bezier{50}(140.00,85.00)(150.00,92.00)(150.00,100.00)
\bezier{50}(150.00,100.00)(150.00,108.00)(140.00,115.00)
\bezier{100}(140.00,115.00)(125.00,125.00)(100.00,125.00)
\bezier{100}(100.00,125.00)(75.00,125.00)(60.00,115.00)
\bezier{50}(60.00,115.00)(50.00,108.00)(50.00,100.00)
\
\
\smalllines
\bezier{300}(67.00,103.00)(101.00,74.00)(133.00,103.00)
\bezier{250}(70.00,100.00)(100.00,126.00)(130.00,100.00)
\
\
\bezier{150}(30.00,100.00)(30.00,122.00)(50.00,130.00)
\bezier{180}(50.00,130.00)(68.00,140.00)(100.00,140.00)
\bezier{150}(30.00,100.00)(30.00,78.00)(50.00,70.00)
\bezier{180}(50.00,70.00)(68.00,60.00)(100.00,60.00)
\bezier{150}(170.00,100.00)(170.00,78.00)(150.00,70.00)
\bezier{180}(150.00,70.00)(132.00,60.00)(100.00,60.00)
\
\
\fatlines
\bezier{50}(280.00,99.00)(280.00,91.00)(290.00,84.00)
\bezier{100}(290.00,84.00)(306.00,74.00)(330.00,74.00)
\bezier{100}(330.00,74.00)(355.00,74.00)(370.00,84.00)
\bezier{50}(370.00,84.00)(380.00,91.00)(380.00,99.00)
\bezier{50}(380.00,99.00)(380.00,107.00)(370.00,114.00)
\bezier{100}(370.00,114.00)(355.00,124.00)(330.00,124.00)
\bezier{100}(330.00,124.00)(305.00,124.00)(290.00,114.00)
\bezier{50}(290.00,114.00)(280.00,107.00)(280.00,99.00)
\
\
\smalllines\bezier{200}(297.00,102.00)(331.00,73.00)(363.00,102.00)
\bezier{200}(300.00,99.00)(330.00,125.00)(360.00,99.00)
\bezier{150}(260.00,99.00)(260.00,77.00)(280.00,69.00)
\bezier{260}(280.00,69.00)(298.00,59.00)(330.00,59.00)
\bezier{120}(400.00,99.00)(400.00,77.00)(380.00,69.00)
\bezier{160}(380.00,69.00)(362.00,59.00)(330.00,59.00)
\
\
\bezier{176}(400.00,99.00)(400.00,121.00)(380.00,129.00)
\bezier{212}(380.00,129.00)(362.00,139.00)(330.00,139.00)
\
\
\
\fatlines
\bezier{300}(150.00,100.00)(150.00,240.00)(215.00,252.00)
\bezier{300}(280.00,100.00)(280.00,240.00)(215.00,252.00)
\
\smalllines
\bezier{400}(198.00,252.00)(144.00,233.00)(143.00,115.00)
\bezier{400}(287.00,115.00)(287.00,233.00)(231.00,252.00)
\bezier{400}(158.00,115.00)(158.00,220.00)(205.00,244.00)
\bezier{400}(273.00,115.00)(275.00,210.00)(225.00,242.00)
\bezier{28}(198.00,252.00)(202.00,252.00)(204.00,250.00)
\bezier{28}(204.00,250.00)(206.00,247.00)(205.00,244.00)
\bezier{28}(231.00,252.00)(227.00,252.00)(226.00,249.00)
\bezier{28}(226.00,249.00)(225.00,246.00)(225.00,242.00)
\bezier{40}(143.00,115.00)(144.00,110.00)(149.00,110.00)
\bezier{52}(149.00,110.00)(158.00,110.00)(158.00,115.00)
\bezier{44}(273.00,115.00)(273.00,110.00)(280.00,110.00)
\bezier{44}(281.00,110.00)(287.00,110.00)(287.00,115.00)
\
\
\bezier{150}(100.00,140.00)(128.00,140.00)(143.00,132.00)
\bezier{108}(170.00,100.00)(170.00,117.00)(158.00,126.00)
\bezier{150}(330.00,139.00)(301.00,139.00)(287.00,132.00)
\bezier{102}(260.00,99.00)(260.00,117.00)(273.00,125.00)
\bezier{92}(192.00,250.00)(193.00,264.00)(200.00,269.00)
\bezier{68}(200.00,269.00)(207.00,275.00)(215.00,275.00)
\bezier{80}(215.00,275.00)(223.00,275.00)(232.00,267.00)
\bezier{84}(232.00,237.00)(226.00,229.00)(215.00,229.00)
\bezier{92}(197.00,239.00)(205.00,229.00)(215.00,229.00)
\bezier{80}(232.00,267.00)(239.00,259.00)(237.00,250.00)
\put(100.00,152.00){\makebox(0,0)[cb]{$c_1$}}
\put(330.00,152.00){\makebox(0,0)[cb]{$c_2$}}
\put(215.00,254.00){\makebox(0,0)[cb]{$p$}}
\put(100.00,150.00){\vector(0,-1){25}}
\put(330.00,150.00){\vector(0,-1){26}}
\end{picture}

\caption{An expressway with two generators}
\end{figure}
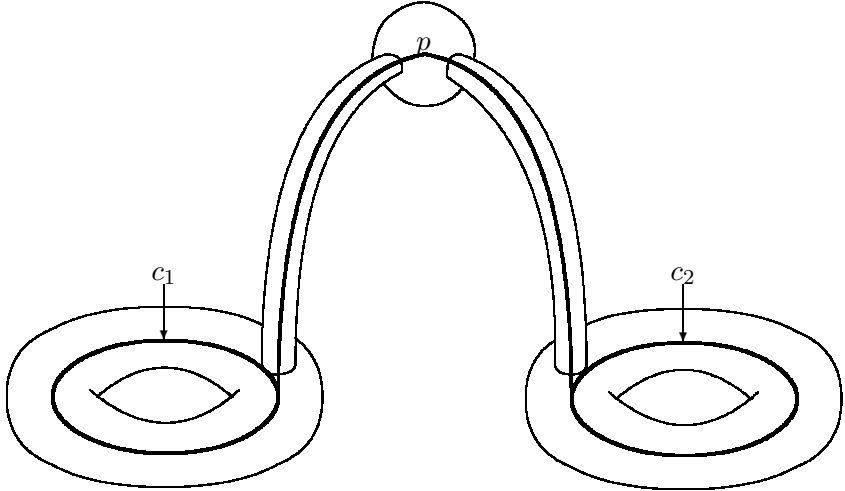

When the author tried to generalise Bangert's construction of
Hedlund metrics (\cite{bangertmin}) to manifolds with arbitrary 
fundamental groups, he took closed curves $c_1,\dots,c_k$ based
in $p\in M$, whose homotopy classes $[c_i]$ generate $\Pi_1$.
By a small perturbation it is possible to 
transform the $c_i$ into smooth disjoint embeddings of 
$S^1\hookrightarrow M$ passing near $p$ if $\dim M\geq 3$.
Now he chose a Riemannian metric that is very small in a small 
neighborhood of the $c_i$ and small in the neighborhood of certain paths 
joining the $c_i$ to $p$, but relatively big outside these neighborhoods.
We can assume that the $c_i$ are hyperbolic minimal geodesics of length $\ep$.

Roughly speaking, the Riemannian distance looks like the distance,
a car driver has in his mind: there are some ``expressways'' 
(the neighborhoods of the $c_i$ and the joining paths) where normed curves
run very fast and in other regions where they move relatively slow.
So minimal geodesics run most of their time on these ``expressways''.
To be more precise the author showed that 
if $c$ is a minimal geodesic and $E$ the expressway, then the length of
each connected component of $c^{-1}(M\ohne E)$ is small.
But we do not know whether the total length of $c^{-1}(M\ohne E)$ is bounded.

So for arbitrary nilpotent fundamental groups the author was 
unable to get analogues of propositions \ref{beschrlaeng} and \ref{beschrtub}.
The situation is much nicer when $\Pi_1$ is of bounded minimal 
generation with respect to the $[c_i]$ (Definition \ref{bmg}).
Here we get analogues of propositions \ref{beschrlaeng}
and \ref{beschrtub}. 
The classification of minimal geodesics then can be reduced to 
combinatorial group properties of $\Pi_1$, and every minimal geodesic 
is asymptotic in each of its senses to one of the $c_i$.
So the results are similar as on $T^n$ ($n\geq 3$).
Unfortunately the bounded-minimal-generation-problem seems to be hard.

\section{Groups of Bounded Minimal Generation}\label{secbmg}

Let $S_\Ga$ be a finite set of generators of the group $\Ga$, i.e.\ 
every $\ga\in \Ga $ can be written as a \me{word} $\wordl sl$
with $s_i\in \SGA$. The number $l=:l(\word{s})$ 
is the \me{length} of the word, furthermore we define the change number
\begin{eqnarray*}
 \CC(\wordl sl)&:=&\# \big\{i \in \{1,2,\dots,l-1\} \, | \, 
    s_i\neq s_{i+1}\big\}+1\\
  &=& \sup \left\{ k\,|\, s_{i_1}\neq \dots \neq s_{i_k},\,\, 
      1\leq i_1<\dots<i_k\leq l\right\}.
\end{eqnarray*}
This is a group theoretical analogue of Definition \ref{changenumber}.
For the empty word representing the neutral element 
we set $l(\emptyset)=\CC(\emptyset)=0$.

The word $\word{s}$ is 
of minimal length if every 
$\word{s'}$ representing the same $\ga\in \Ga$ satisfies:
  $$l(\word{s})\leq l(\word{s'}).$$

\begin{definition}\label{bmg}
$(\Ga,S_\Ga)$ is a \me{Group of Bounded Minimal Generation (BMG group)}
if there is a $B\in \natur$ such that 
every $\ga \in \Ga$ can be represented by a word of minimal length
$\word{s}$ in $\SGA$ with $\CC(\word{s})\leq B$. The minimal such $B$ 
will be called the \me{bound}.
\end{definition}

Every finitely generated abelian group together with an arbitrary 
finite set of generators $S_\Ga$ is a BMG group with $B\leq \#S_\Ga$.

Gromov proved in \cite{gromov1} that every finitely generated group of 
polynomial growth is virtually nilpotent, i.e.\ it contains 
a nilpotent subgroup of finite index. If $(\Ga,S_\Ga)$ is a BMG group,
then $\Ga$ is of polynomial growth and therefore virtually nilpotent.
Yet, it is not clear to the author if the converse holds.

\begin{openproblem}
Is every finitely generated virtually nilpotent group a BMG group?
\end{openproblem}
\remark
It is not even clear whether the BMG property is independent of the 
choice of the set of generators.
\bigskip

\head{Constructing new BMG groups from old ones.}
It is staightforward to show:
\begin{enumerate}
\item If $(\Ga,S_\Ga)$ and $(\Ga',S_{\Ga'})$ are BMG groups with bounds $B$ 
and $B'$, 
then $(\Ga\times \Ga', S_\Ga \cup S_{\Ga'})$ is a BMG group with bound 
$\leq B+B'$.
\item If $(\Ga, S_\Ga)$ is a BMG-group with bound $B$ 
and $h\colon \Ga \to \Ga'$ a 
group homomorphism, then $(h(\Ga),h(S_\Ga))$ is a 
BMG-group with bound $\leq B$.
\item If $\Ga_2$ is the semidirect product of a BMG group $(\Ga_1,S_{\Ga_1})$ 
and a finite group,
then there is a generating system 
$S_{\Gamma_2}$ of $\Gamma_2$ such that $(\Gamma_2,S_{\Gamma_2})$ 
is a BMG group.
\end{enumerate}

\noindent
\head{Heisenberg groups.} 
For $m \in \natur$, $\vp,\vq \in \reell^m$, $z \in \reell$ we define
the matrix
  $$M(\vp,\vq,z):=\pmatrix{1 & \vp^{\,\,t} & z \cr 0 & 1_m & \vq \cr 0 & 0 & 1}$$
$\hei{m}:=\big\{ M(\vp,\vq,z)\, |\, \vp,\vq 
 \in \reell^m,\ z \in \reell\big\}$
is the $2m+1$-dimensional Heisenberg group.

For $r=(r_1,\dots,r_m)\in \natur^m$ such that $r_i$ divides $r_{i+1}$,
$1\leq i < m$, we set 
  $$r \Z^m := \left\{ (r_1 x_1,\dots,r_mx_m)\, | \, x_j \in \Z \right\}$$ 
  $$\Ga_r :=\left\{ M(\vp,\vq,z) \, |\, \vp \in r \Z^m, 
    \vq \in \Z^m, z \in \Z\right\}.$$
These $\Ga_r$ are \me{lattices} 
in $\hei{m}$, i.e.\ discrete, cocompact subgroups.
 
\begin{theorem}[\cite{gowi}, \S 2.]
For every lattice $\Ga$ of $\hei{m}$ there exists
a unique $r$ and an automorphism of $\hei{m}$ mapping $\Ga$ to $\Ga_r$.
\end{theorem}    

\begin{theorem}\label{allebmg}
Every lattice $\Ga$ of $\hei{m}$
has a set of generators $S$, such that $(\Ga,S)$ is a BMG-group. 
\end{theorem}

\remark 
It is not difficult to show that any discrete subgroup of $\hei{m}$
is of the form
$\Gamma' \times \Z^ k$, where $\Gamma'$ is trivial or a lattice in a 
Heisenberg group. So Theorem \ref{allebmg} immediately generalizes to
discrete subgroups. 
\bigskip

\begin{lemma}\label{einebmg}
$\Ga_1= \big\{ M(p,q,z)\, |\, p,q,z \in \Z\big\}\subset \HH_1$ together with
\newline $S_{\Ga_1}:=\left\{M(1,0,0), M(0,1,0)\right\}$ 
is a BMG-group.
\end{lemma}

\spacebetweentheoremandproof

\begin{pf*}{Proof of Theorem \ref{allebmg}}
We can assume $\Ga=\Ga_r$. We denote the standard basis of $\reell^m$ by
$e_1,\dots,e_m$.
The mappings 
\begin{eqnarray*}
  f_i\colon \hei{1} & \to & \hei{m} \cr
  M(p,q,z)          & \mapsto & M(pr_ie_i,qe_i,zr_i)
\end{eqnarray*}
for $i=1,\dots , m$ and 
\begin{eqnarray*}
  f_0 \colon \reell & \to & \hei{m} \cr
              z & \mapsto & M(0,0,z)   
\end{eqnarray*}
define a group epimorphism 
\begin{eqnarray*}
  f \colon \reell \times \bigtimes_{i=1}^m \hei{1} & \to & \hei{m}\cr
  (z,h_1,\dots,h_m) & \mapsto & f_0(z)f_1(h_1)f_2(h_2)\dots f_m(h_m)
\end{eqnarray*}
that maps $\Z \times \bigtimes_{i=1}^m \Ga_1$ to $\Ga_r$.
Using lemma \ref{einebmg} we know that $\left(\Ga_1,S_{\Ga_1}\right)$ 
is a BMG-group. Applying the above constructions
of new BMG-groups we see that $(\Ga,S)$ is a BMG-group 
for $S:=f(\bigcup_{i=0}^m S_i)$, where $S_0$ is $1\in \reell$ 
and $S_i$ is $S_{\Ga_1}$ in the $i$-th slot.
\qed
\end{pf*}

\begin{pf*}{Proof of Lemma \ref{einebmg}}
We set $g_1:=M(1,0,0)$, $g_2:=M(0,1,0)$, then 
$[g_1,g_2]=g_1g_2g_1^{-1}g_2^{-1}=M(0,0,1)$.
For every word $w$ in the generators $g_1$ and $g_2$ 
we will explain how to construct 
a word $w'$ of minimal length and with $\CC(w')\leq 6$ representing the same 
$\ga \in \Ga$.

In order to construct $w'$ we will give a geometric 
interpretation of the problem.
To every word $w$ we will associate inductively a path $p(w)$ in 
$$P:=\left\{\zvec xy\in \reell^2 \, | \, x \in \Z \mbox{ or } 
 y \in \Z\right\}.$$
At first we associate to $w=\empt$ the constant path in $\zvec 00$, and 
to $w=g_1$ (resp.\ $w=g_2,g_1^{-1},g_2^{-1}$) the straight line from 
$\zvec 00$ to $\zvec 10$ (resp.\ $\zvec 01,
\linebreak[0] \zvec{-1}0,\linebreak[0] \zvec0{-1}$).
Then we associate to the word $w_1w_2$ the path $p(w_1w_2)$
that consists of the path $p(w_1)$ and then $p(w_2)$, translated 
by the endpoint of $p(w_1)$ --- 
we get a path from $\zvec 00$ to the sum 
of the endpoints of $p(w_1)$ and $p(w_2)$.

Now we define 
  $$I(w)=\int_{p(w)} x\, dy.$$ 
Let $(i_1,i_2)$ be the endpoint of $p_w$.
Then $w$ represents ${g_2}^{i_2}{g_1}^{i_1} \left[g_1,g_2\right]^{I(w)}$,
$l(w)=\length (p(w))$, the change number $\CC(w)$ of $w$ is equal to the 
number of direction changes of $p(w)$ plus 1. 

We get an geometric interpretation of $I(w)$ by applying Stoke's theorem:
  $$I(w)=\int_{p(w)} x \, dy = \int _{p\left(w{g_1}^{-i_1}{g_2}^{-i_2}\right)}
    x \, dy = \int_A dx \wedge dy,$$
where $A$ is the 2-chain whose boundary 
is $p\left(w{g_1}^{-i_1}{g_2}^{-i_2}\right)$. 
As $dx \wedge dy$ is the volume element of $\reell^2$, we interpret
$I(w)$ as the oriented area of $A$.

We get an isoperimetric problem on $P$: Finding words of minimal length
means finding paths of minimal length in $P$ with given integral $\int x\,dy$.

With this geometric interpretation the lemma is almost obvious.
Using symmetries we can assume $i_1,i_2\geq 0$, $I(w)\geq 0$.
Consider the special cases
\begin{enumerate}
\item $I(w)\leq i_1\cdot i_2$,
\item $i_1i_2<I(w)\leq \max\left\{i_1,i_2\right\}^2$,
\item $\max\left\{i_1,i_2\right\}^2 < I(w)$.
\end{enumerate}
In each case we can find a word $w'$ of minimal length, equivalent to $w$ 
with $\CC(w')\leq 6$.
\qed
\end{pf*}
\remark
Michael Stoll \cite{stoll} treats the continuous 
analogue of bounded minimal generation for nilpotent Lie-groups.
He proves that every 2-step nilpotent Lie-group fulfills it, 
but he states that there are counterexamples 
for 3-step nilpotent Lie-groups.
\bigskip

\remark
Several authors (\cite{carkel},\cite{tavgen},\cite{rapinchuk1} 
and \cite{rapinchuk2}) consider definitions
(``groups of bounded generation'', ``groups of finite width'') which are 
formally related to our groups of bounded minimal generation.
The main difference is that these notions use arbitrary not only minimal
representatives. 
\bigskip


\end{document}